\documentclass[11pt]{amsart}
\usepackage{amssymb,amsmath,a4wide}
\usepackage{graphicx}
\usepackage[all]{xy}
\usepackage{stmaryrd}

\title[Combinatorial and categorical semantics of process
algebra]{Homotopical equivalence of combinatorial and categorical
  semantics of process algebra} 
\author[P. Gaucher]{Philippe Gaucher}
\address{Laboratoire PPS  (CNRS UMR 7126)\\ Universit{\'e} Paris 7--Denis Diderot\\
  Case 7014\\ 75205 PARIS Cedex 13 \\ France}
\email{gaucher@pps.jussieu.fr}
\urladdr{http://www.pps.jussieu.fr/{\~{}}gaucher/} 
\subjclass{55U35,18G55,68Q85} 
\keywords{homotopy colimit, homotopy limit, weak colimit, weak limit, precubical set, time flow,  higher dimensional automaton, process algebra, concurrency}

\swapnumbers

\newcommand{\C}{\mathcal{C}}

\newcommand{\N}{\mathbb{N}}

\newcommand{\de}{\partial}
\newcommand{\p}\times
\renewcommand{\vec}{\overrightarrow}
\renewcommand{\P}{\mathbb{P}}

\newcommand{\be}{\begin{equation}}
\newcommand{\ee}{\end{equation}}
\newcommand{\bea}{\begin{eqnarray}}
\newcommand{\eea}{\end{eqnarray}}
\newcommand{\beas}{\begin{eqnarray*}}
\newcommand{\eeas}{\end{eqnarray*}}

\newtheorem{thm}{Theorem}[section]
\newtheorem{prop}[thm]{Proposition}

\newtheorem{cor}[thm]{Corollary}

\newtheorem{defn}[thm]{Definition}
\newtheorem{nota}[thm]{Notation}
\newtheorem{rem}[thm]{Remark}
\newcommand{\bd}{\begin{defn}}
\newcommand{\ed}{\end{defn}}
\newcommand{\bp}{\begin{prop}}
\newcommand{\ep}{\end{prop}}
\newcommand{\bth}{\begin{thm}}
\renewcommand{\eth}{\end{thm}}
\newcommand{\bpf}{\begin{proof}}
\newcommand{\epf}{\end{proof}}

\newcommand{\fl}[1]{\ar@{->}[ll]_-{#1}}
\newcommand{\fr}[1]{\ar@{->}[rr]^-{#1}}
\newcommand{\fd}[1]{\ar@{->}[dd]_-{#1}}
\newcommand{\fu}[1]{\ar@{->}[uu]^-{#1}}
\newcommand{\f}[2]{\ar@{->}[#1]|{#2}}
\newcommand{\ff}[2]{\ar@2{->}[#1]|{#2}}
\newcommand{\frr}[1]{\ar@{->}[rrrr]^-{#1}}

\renewcommand{\top}{{\mathbf{Top}}}
\newcommand{\ho}{{\mathbf{Ho}}}
\newcommand{\iso}{\cong}

\newcommand{\ot}{\otimes}

\renewcommand{\leq}{\leqslant}
\renewcommand{\geq}{\geqslant}

\def\cartesien{%
  \ar@{-}[]+R+<6pt,-2pt>;[]+RD+<6pt,-6pt>%
  \ar@{-}[]+D+<2pt,-6pt>;[]+RD+<6pt,-6pt>%
}
\def\cocartesien{%
  \ar@{-}[]+L+<-6pt,+2pt>;[]+LU+<-6pt,+6pt>%
  \ar@{-}[]+U+<-2pt,+6pt>;[]+LU+<-6pt,+6pt>%
}
\def\hocartesien{%
  \ar@{-}[]+R+<6pt,-2pt>;[]+RD+<6pt,-6pt>_{h}%
  \ar@{-}[]+D+<2pt,-6pt>;[]+RD+<6pt,-6pt>%
}
\def\hococartesien{%
  \ar@{-}[]+L+<-6pt,+2pt>;[]+LU+<-6pt,+6pt>_{h}%
  \ar@{-}[]+U+<-2pt,+6pt>;[]+LU+<-6pt,+6pt>%
}
\def\wcartesien{%
  \ar@{-}[]+R+<6pt,-2pt>;[]+RD+<6pt,-6pt>_{w}%
  \ar@{-}[]+D+<2pt,-6pt>;[]+RD+<6pt,-6pt>%
}
\def\wcocartesien{%
  \ar@{-}[]+L+<-6pt,+2pt>;[]+LU+<-6pt,+6pt>_{w}%
  \ar@{-}[]+U+<-2pt,+6pt>;[]+LU+<-6pt,+6pt>%
}

\newcommand{\brm}[1]{\rm{\mathbf{#1}}}

\newcommand{\dtop}{{\brm{Flow}}}
\newcommand{\set}{{\brm{Set}}}
\newcommand{\poset}{{\brm{PoSet}}}
\newcommand{\tdtop}{{\brm{FLOW}}}

\newcommand{\proc}{{\brm{Proc}}}
\newcommand{\glob}{{\rm{Glob}}}
\DeclareMathOperator{\map}{Map}
\DeclareMathOperator{\rec}{rec}
\DeclareMathOperator{\id}{Id}

\DeclareMathOperator{\Obj}{Obj}
\newcommand{\liminj}{\varinjlim}
\newcommand{\limproj}{\varprojlim}

\makeatletter
\def\varholim@#1#2{%
  \vtop{\m@th\ialign{##\cr
    \hfil$#1\operator@font holim$\hfil\cr
    \noalign{\nointerlineskip\kern1.5\ex@}#2\cr
    \noalign{\nointerlineskip\kern-\ex@}\cr}}%
}
\def\holimproj{%
  \mathop{\mathpalette\varholim@{\leftarrowfill@\textstyle}}\nmlimits@
}
\def\holiminj{%
  \mathop{\mathpalette\varholim@{\rightarrowfill@\textstyle}}\nmlimits@
}
\def\varwlim@#1#2{%
  \vtop{\m@th\ialign{##\cr
    \hfil$#1\operator@font wlim$\hfil\cr
    \noalign{\nointerlineskip\kern1.5\ex@}#2\cr
    \noalign{\nointerlineskip\kern-\ex@}\cr}}%
}
\def\wlimproj{%
  \mathop{\mathpalette\varwlim@{\leftarrowfill@\textstyle}}\nmlimits@
}
\def\wliminj{%
  \mathop{\mathpalette\varwlim@{\rightarrowfill@\textstyle}}\nmlimits@
}
\makeatother

\DeclareMathOperator{\COSK}{COSK}

\newcommand{\sis}{\Delta^{op}\set}
\DeclareMathOperator{\cell}{{\brm{cell}}}
\DeclareMathOperator{\cof}{{\brm{cof}}}
\DeclareMathOperator{\inj}{{\brm{inj}}}
\newcommand{\ddownarrow}{{\downarrow}}
\DeclareMathOperator{\dgm}{dgm}
\begin{document}

\begin{abstract}
  It is possible to translate a modified version of K. Worytkiewicz's
  combinatorial semantics of CCS (Milner's Calculus of Communicating
  Systems) in terms of labelled precubical sets into a categorical
  semantics of CCS in terms of labelled flows using a geometric
  realization functor. It turns out that a satisfactory semantics in
  terms of flows requires to work directly in their homotopy category
  since such a semantics requires non-canonical choices for
  constructing cofibrant replacements, homotopy limits and homotopy
  colimits. No geometric information is lost since two precubical sets
  are isomorphic if and only if the associated flows are weakly
  equivalent. The interest of the categorical semantics is that
  combinatorics totally disappears.  Last but not least, a part of the
  categorical semantics of CCS goes down to a pure homotopical
  semantics of CCS using A. Heller's privileged weak limits and
  colimits. These results can be easily adapted to any other process
  algebra for any synchronization algebra.
\end{abstract}

\maketitle

\tableofcontents

\section{Introduction}

This paper is the companion paper of \cite{ccsprecub}. The preceding
paper was devoted to fixing K. Worytkiewicz's combinatorial semantics
of CCS (Milner's Calculus of Communicating System) \cite{0683.68008}
\cite{MR1365754} in terms of labelled precubical sets \cite{exHDA} in
order to stick to the higher dimensional automata paradigm. This
paradigm states that the concurrent execution of $n$ actions must be
abstracted by exactly \textit{one} full $n$-cube: see \cite{ccsprecub}
Theorem~5.2 for a rigorous formalization of this paradigm and also
Proposition~\ref{pourquoi_non_degenere} of this paper. There was a
problem in K.  Worytkiewicz's approach because of a version of the
labelled coskeleton construction adding too many cubes and therefore
not satisfying Proposition~\ref{pourquoi_non_degenere}.  The purpose
of the preceding paper was also to built an appropriate geometric
realization functor from labelled precubical sets to labelled flows.
The little bit surprising fact arising from this construction was that
a satisfactory geometric realization functor does require the use of
the model structure of flows introduced in \cite{model3}. A
consequence of the preceding paper was to give a proof of the
expressiveness of the category of flows.  The geometric intuition
underlying these two semantics, i.e.  in terms of precubical sets and
in terms of flows, is extensively explained in Section~\ref{introbis}
which must be considered as a part of this introduction.

In this work, we push a little bit further the study of the semantics
of CCS in terms of labelled flows. Indeed, we explain the effect of
the geometric realization functor on each operator of CCS.
Section~\ref{effect} is the section of the paper presenting these new
results. In particular, Theorem~\ref{parho} presents an interpretation
of the parallel composition with synchronization in terms of flows
without any combinatorial construction. This must be considered as the
main result of the paper.

The only case treated in this paper is the one of CCS without message
passing. But all the results can be easily adapted to any other
process algebra with any other synchronization algebra. The case of
TCSP \cite{0628.68025} was explained in \cite{ccsprecub}. For general
synchronization algebras, all proofs of the paper are exactly the
same, except the proof of Proposition~\ref{p1eq} which must be very
slightly modified: see the comment in the footnote~\ref{moregeneral}.

\subsection*{Outline of the paper}

Section~\ref{introbis} explains, with the example of the concurrent
execution of two actions $a$ and $b$, the geometric intuition
underlying the two semantics studied in this paper. It must be
considered as part of the introduction and it is strongly recommended
the reading for anyone not knowing the subject (and also for the other
ones). In particular, the notions of labelled precubical set and of
labelled flow are reminded here. Section~\ref{comb} recalls the syntax
of CCS and the construction of the combinatorial semantics of
\cite{ccsprecub} in terms of labelled precubical sets. The geometric
realization functor is then introduced in Section~\ref{catsem}. Since
we \textit{do} need to work in the homotopy category of flows,
Section~\ref{relevance} proving that two precubical sets are
isomorphic if and only if the associated flows are weakly S-homotopy
equivalent is fundamental.  Finally Section~\ref{effect} is an
exposition of the effect of the geometric realization functor from
precubical sets to flows on each operator defining the syntax of CCS.
It is the technical core of the paper.  And Section~\ref{bonus} is a
bonus explaining some ideas towards a pure homotopical semantics of
CCS: Theorem~\ref{weakweakweak} is a consequence of all the theorems
of Section~\ref{effect} and of some known facts about realization of
homotopy commutative diagrams over free Reedy categories and their
links with some kinds of weak limits and weak colimits in the homotopy
category of a model category.

\subsection*{Prerequisites}

The reading of this work requires some familiarity with model category
techniques \cite{MR99h:55031} \cite{ref_model2}, with category theory
\cite{MR1712872} \cite{MR96g:18001a}\cite{gz}, and also with locally
presentable categories \cite{MR95j:18001}. We use the locally
presentable category of $\Delta$-generated topological spaces.
Introductions about these spaces are available in \cite{delta}
\cite{FR} and \cite{interpretation-glob}.

\subsection*{Notations}

Let $\C$ be a cocomplete category. The class of morphisms of $\C$ that
are transfinite compositions of pushouts of elements of a set of
morphisms $K$ is denoted by $\cell(K)$. An element of $\cell(K)$ is
called a \textit{relative $K$-cell complex}.  The category of sets is
denoted by $\set$. The class of maps satisfying the right lifting
property with respect to the maps of $K$ is denoted by $\inj(K)$. The
class of maps satisfying the left lifting property with respect to the
maps of $\inj(K)$ is denoted by $\cof(K)$. The cofibrant replacement
functor of a model category is denoted by $(-)^{cof}$. The notation
$\simeq$ means \textit{weak equivalence} or \textit{equivalence of
  categories}, the notation $\iso$ means \textit{isomorphism}. The
notation $\id_A$ means identity of $A$.  The initial object (resp.
final object) of a category is denoted by $\varnothing$ (resp.
$\mathbf{1}$). The cofibrant replacement functor of any model category
is denoted by $(-)^{cof}$.  The category of partially ordered set or
poset together with the strictly increasing maps ($x<y$ implies
$f(x)<f(y)$) is denoted by $\poset$. The set of morphisms from an
object $X$ to an object $Y$ of a category $\C$ is denoted by
$\C(X,Y)$.

\subsection*{Acknowledgments} I thank very much Andrei
R\u{a}dulescu-Banu for bringing Theorem~\ref{weak} to my attention.

\section{Example of two concurrent executions}
\label{introbis}

We want to explain in this section what we mean by
\textit{combinatorial semantics} and \textit{categorical semantics}
with the example of the concurrent execution of two actions $a$ and
$b$. This section also recalls the definitions of \textit{labelled
  precubical set} and of \textit{labelled flow}. For other references
about topological models of concurrency, see \cite{survol} for a
survey.

\begin{figure}
\begin{center}
\includegraphics[width=5cm]{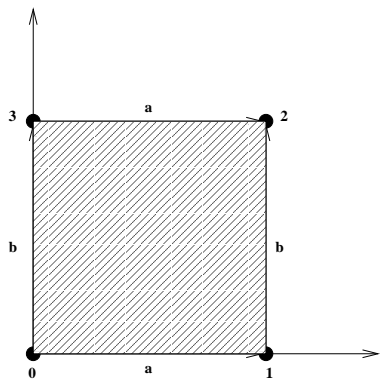}
\end{center}
\caption{Concurrent execution of two actions $a$ and $b$} 
\label{carreplein}
\end{figure}

Consider two actions $a$ and $b$ whose concurrent execution is
topologically represented by the square $[0,1]^2$ of
Figure~\ref{carreplein}.  The topological space $[0,1]^2$ itself
represents the underlying state space of the process. Four
distinguished states are depicted on Figure~\ref{carreplein}. The
state $0 = (0,0)$ is the initial state. The state $2 = (1,1)$ is the
final state.  At the state $1 = (1,0)$, the action $a$ is finished and
the action $b$ is not yet started.  At the state $3 = (0,1)$, the
action $b$ is finished and the action $a$ is not yet started. So the
boundary $[0,1]\p\{0,1\} \cup \{0,1\}\p [0,1]$ of the square $[0,1]^2$
models the sequential execution of the actions $a$ and $b$ whereas
their concurrent execution is modeled by including the $2$-dimensional
square $]0,1[ \p ]0,1[$. In fact, the possible execution paths from
the initial state $0 = (0,0)$ to the final state $2 = (1,1)$ are all
continuous paths from $0 = (0,0)$ to $2 = (1,1)$ which are
non-decreasing with respect to each axis of coordinates.
Nondecreasingness corresponds to irreversibility of time.

In the combinatorial semantics, the preceding situation is abstracted
by a $2$-cube viewed as a precubical set. Let us now recall the
definition of these objects. A good reference for presheaves is
\cite{MR1300636}.

\begin{nota}
  Let $[0] = \{0\}$ and $[n] = \{0,1\}^n$ for $n \geq 1$.  By
  convention, $\{0,1\}^0=\{0\}$.
\end{nota}

Let $\delta_i^\alpha : [n-1] \rightarrow [n]$ be the set map defined
for $1\leq i\leq n$ and $\alpha \in \{0,1\}$ by
$\delta_i^\alpha(\epsilon_1, \dots, \epsilon_{n-1}) = (\epsilon_1,
\dots, \epsilon_{i-1}, \alpha, \epsilon_i, \dots, \epsilon_{n-1})$.
The small category $\square$ is by definition the subcategory of the
category of sets with set of objects $\{[n],n\geq 0\}$ and generated
by the morphisms $\delta_i^\alpha$.

\bd \cite{Brown_cube} The category of presheaves over $\square$,
denoted by $\square^{op}\set$, is called the category of {\rm
  precubical sets}.  A precubical set $K$ consists of a family of sets
$(K_n)_{n \geq 0}$ and of set maps $\de_i^\alpha:K_n \rightarrow
K_{n-1}$ with $1\leq i \leq n$ and $\alpha\in\{0,1\}$ satisfying the
cubical relations $\de_i^\alpha\de_j^\beta = \de_{j-1}^\beta
\de_i^\alpha$ for any $\alpha,\beta\in \{0,1\}$ and for $i<j$. An
element of $K_n$ is called a {\rm $n$-cube}. \ed

Let $\square[n]:=\square(-,[n])$. By the Yoneda lemma, one has the
natural bijection of sets \[\square^{op}\set(\square[n],K)\iso K_n\]
for every precubical set $K$. The boundary of $\square[n]$ is the
precubical set denoted by $\de \square[n]$ defined by removing the
interior of $\square[n]$:
\begin{itemize} 
\item $(\de \square[n])_k := (\square[n])_k$ for $k<n$
\item $(\de \square[n])_k = \varnothing$ for $k\geq n$.
\end{itemize} 
In particular, one has $\de \square[0] = \varnothing$.

So the $2$-cube $\square[2]$ models the underlying time flow of the
concurrent execution of $a$ and $b$. However, the same $2$-cube models
the underlying time flow of the concurrent execution of any pair of
actions. So we need a notion of labelling.

Let $\Sigma$ be a set of labels, containing among other things the two
actions $a$ and $b$.

\bp \cite{labelled} \label{cubeetiquette} Put a total ordering $\leq$
on $\Sigma$.  Let
\begin{itemize}
\item $(!\Sigma)_0=\{()\}$ (the empty word)
\item for $n\geq 1$, $(!\Sigma)_n=\{(a_1,\dots,a_n)\in \Sigma \p \dots
  \p \Sigma,a_1 \leq \dots \leq a_n \}$
\item $\de_i^0(a_1,\dots,a_n) = \de_i^1(a_1,\dots,a_n) =
  (a_1,\dots,\widehat{a_i},\dots,a_n)$ where the notation
  $\widehat{a_i}$ means that $a_i$ is removed.
\end{itemize}
Then these data generate a precubical set. \ep 

\begin{rem} The isomorphism class of $!\Sigma$ does not depend of the
  choice of the total ordering on $\Sigma$. \end{rem}

\bd (Goubault) A {\rm labelled precubical set} is an object of the
comma category $\square^{op}\set \ddownarrow !\Sigma$. \ed

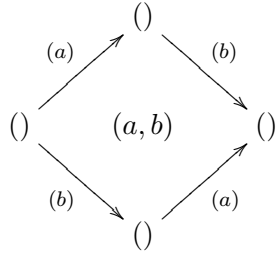
\begin{figure}
\[
\xymatrix{
& () \ar@{->}[rd]^{(b)}&\\
()\ar@{->}[ru]^{(a)}\ar@{->}[rd]_{(b)} & (a,b) & ()\\
&()\ar@{->}[ru]_{(a)}&}
\]
\caption{Concurrent execution of $a$ and $b$ with  $a\leq b$ as labelled precubical set}
\label{exab}
\end{figure}

In the combinatorial semantics, the concurrent action of the two
actions $a$ and $b$ is then modeled by the labelled $2$-cube
$\ell:\square[2]\rightarrow !\Sigma$ sending the identity of $[2]$
(the interior of the square) to $(a,b)$ if $a\leq b$ or to $(b,a)$ if
$b\leq a$ as depicted in Figure~\ref{exab}.

The categorical semantics is much simpler to explain. Each of the four
distinguished states $0$, $1$, $2$ and $3$ of Figure~\ref{carreplein}
is represented by an object of a small category. Each execution path
of the boundary of the square is represented by a morphism with the
composition of execution paths corresponding to the composition of
morphisms. So one has $6$ execution paths $\vec{01}$, $\vec{12}$,
$\vec{012}$, $\vec{03}$, $\vec{32}$ and $\vec{032}$ with the algebraic
rules $\vec{01} * \vec{12} = \vec{012}$ and $\vec{03} * \vec{32} =
\vec{032}$ where $*$ is of course the composition law. The interior of
the square is then modeled by the algebraic relation $\vec{012} =
\vec{032}$. This small category is nothing else but the small category
corresponding to the poset $\{\widehat{0} < \widehat{1}\}^2$. And this
poset is nothing else but the poset of vertices of the $2$-cube. The
partial ordering models observable time ordering.

In fact, one needs to work with categories \textit{enriched over
  topological spaces} in the sense of \cite{MR2177301}, i.e. with
topologized homsets, for being able to model more complicated
situations of concurrency.  In the situation above, the space of
morphisms is of course discrete.  For various mathematical reasons,
e.g. \cite{model3} Section~20 and \cite{diCW} Section~6, one also
needs to work with small category \textit{without identity maps}.
Note that the category of small categories without identity maps and
the usual one of small categories are certainly not equivalent since a
round-trip using the adjunction between them adds a loop to each
object. And a lots of theorems proved in the framework of flows (i.e.
small categories without identity maps enriched over topological
spaces) are merely wrong whenever identity maps are added.

In this paper, one will also work with the locally presentable
category of \textit{$\Delta$-generated topological spaces}, denoted by
$\top$, i.e. of spaces which are colimits of simplices.  Several
introductions about these topological spaces are available in
\cite{delta} \cite{FR} and \cite{interpretation-glob} respectively.
Let us only mention one striking property of $\Delta$-generated
topological spaces: as the simplicial sets, they are isomorphic to the
disjoint sum of their [path-]connected components by
\cite{interpretation-glob} Proposition~2.8. This property has a lots
of very nice consequences.

\bd \cite{model3} A {\rm flow} $X$ is a small category without
identity maps enriched over $\Delta$-generated topological spaces.
The composition law of a flow is denoted by $*$.  The set of objects
is denoted by $X^0$.  The space of morphisms from $\alpha$ to $\beta$
is denoted by $\P_{\alpha,\beta} X$~\footnote{Sometimes, an object of
  a flow is called a state and a morphism a (non-constant) execution
  path.}. Let $\P X$ be the disjoint sum of the spaces
$\P_{\alpha,\beta} X$. A morphism of flows $f:X \rightarrow Y$ is a
set map $f^0:X^0 \rightarrow Y^0$ together with a continuous map $\P
f:\P X \rightarrow \P Y$ preserving the structure.  The corresponding
category is denoted by $\dtop$. \ed

Each poset $P$ can be associated with a flow denoted in the same way.
The set of objects $P^0$ is the underlying set of $P$ and there is one
and only one morphism from $\alpha$ to $\beta$ if and only if $\alpha
< \beta$. The composition law is then defined by
$(\alpha,\beta)*(\beta,\gamma) = (\alpha,\gamma)$ for any
$\alpha<\beta<\gamma\in P$. Note that the flow associated with a poset
is \textit{loopless}, i.e. for every $\alpha\in P^0$, one has
$\P_{\alpha,\alpha}P=\varnothing$.  This construction induces a
functor $\poset \rightarrow \dtop$ from the category of posets
together with the strictly increasing maps to the category of flows.

In the categorical semantics, the underlying time flow of the
concurrent execution of two actions $a$ and $b$ is then modeled by the
flow associated with the poset $\{\widehat{0} < \widehat{1}\}^2$.
Like in the combinatorial semantics, one needs a notion of labelling.

\bd \label{defflowlabel} The {\rm flow of labels} $?\Sigma$ is defined
as follows: $(?\Sigma)^0 = \{0\}$ and $\P ?\Sigma$ is the discrete
free associative monoid without unit generated by the elements of
$\Sigma$ and by the algebraic relations $a*b=b*a$ for all $a,b\in
\Sigma$.  \ed

\bd A {\rm labelled flow} is an object of the comma category $\dtop
\ddownarrow ?\Sigma$. \ed

\begin{figure}
\[
\xymatrix{
& 0 \ar@{->}[rd]^{b}&\\
0\ar@{->}[ru]^{a}\ar@{->}[rd]_{b} & = &0 \\
&\ar@{->}[ru]_{a}0&}
\]
\caption{Concurrent execution of $a$ and $b$ as labelled flow}
\label{exab2}
\end{figure}
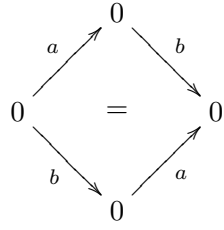

So the concurrent execution of the two actions $a$ and $b$ will be
modeled by the labelled flow $\ell:\{\widehat{0} < \widehat{1}\}^2
\rightarrow ?\Sigma$ with $\ell(\vec{012}) = \ell(\vec{032}) = a*b$,
$\ell(\vec{01}) = a$, $\ell(\vec{12}) = b$, $\ell(\vec{03}) = b$ and
$\ell(\vec{32}) = a$ as in Figure~\ref{exab2}.

\section{Combinatorial semantics of CCS}
\label{comb}

\subsection*{Syntax of the process algebra CCS}

A short introduction about process algebra can be found in
\cite{MR1365754}. An introduction about CCS for mathematician is
available in the companion paper \cite{ccsprecub}. Let $\Sigma$ be a
non-empty set. Its elements are called \textit{labels} or
\textit{actions} or \textit{events}. From now on, the set $\Sigma$ is
supposed to be equipped with an involution $a\mapsto \overline{a}$.
Moreover, the set $\Sigma$ contains a distinct action $\tau$ with
$\tau=\overline{\tau}$. The \textit{process names} are generated by
the following syntax:
\[
P::=nil \ |\ a.P \ |\ (\nu a)P \ |\ P + P \ |\ P|| P \ |\
\rec(x)P(x)\] where $P(x)$ means a process name with one free variable
$x$. The variable $x$ must be \textit{guarded}, that is it must lie in
a prefix term $a.x$ for some $a\in\Sigma$.

\subsection*{Parallel composition with synchronization of labelled
  precubical sets}

\bd \label{def_preshell} Let $\ell:K\rightarrow !\Sigma$ be a labelled
precubical set. Let $n\geq 1$. A {\rm labelled $n$-shell} of $K$ is a
commutative diagram of precubical sets
\[
\xymatrix{
\de\square[n+1] \ar@{->}[rr]^-{x} \fd{} && K \fd{\ell} \\
&& \\
\square[n+1] \fr{} && !\Sigma.}
\]
Suppose moreover that $K_0=[p]$ for some $p\geq 2$. The labelled
$n$-shell above is {\rm non-twisted} if the set map $x_0:[n+1] =
\de\square[n+1]_0 \rightarrow [p] = K_0$ is a composite
\[x_0:[n+1] \stackrel{\phi}\longrightarrow [q]
\stackrel{\psi}\longrightarrow [p]\] where $\psi$ is a morphism of the
small category $\square$ and where $\phi$ is of the form
$(\epsilon_1,\dots,\epsilon_{n+1}) \mapsto
(\epsilon_{i_1},\dots,\epsilon_{i_q})$ such that $1=i_1 \leq \dots
\leq i_q = n+1$ and $\{1,\dots,n+1\}\subset \{i_1,\dots,i_q\}$.
\ed

The map $\phi$ is not necessarily a morphism of the small category
$\square$. For example $\phi:[3] \rightarrow [5]$ defined by
$\phi(\epsilon_1,\epsilon_2,\epsilon_3)=(\epsilon_1, \epsilon_1,
\epsilon_2, \epsilon_3, \epsilon_3)$ is not a morphism of $\square$.
Note the set map $x_0$ is always one-to-one.

All labelled shells of this paper will be supposed non-twisted.

\bd Let $\square_n\subset \square$ be the full subcategory of
$\square$ whose set of objects is $\{[k],k\leq n\}$. The category of
presheaves over $\square_n$ is denoted by $\square_n^{op}\set$. Its
objects are called the {\rm $n$-dimensional precubical sets}.  \ed

Let $K$ be an object of $\square_{1}^{op}\set \ddownarrow !\Sigma$
such that $K_0=[p]$ for some $p\geq 2$.  Let $K^{(n)}$ be the object
of $\square_{n}^{op}\set \ddownarrow !\Sigma$ inductively defined for
$n\geq 1$ by $K^{(1)} = K$ and by the following pushout diagram of
labelled precubical sets (the shells are always non-twisted by
hypothesis):
\[
\xymatrix{
\bigsqcup\limits_{\hbox{labelled $n$-shells}}\de\square[n+1]  \fr{} \fd{}&& K^{(n)} \fd{}\\
&& \\
\bigsqcup\limits_{\hbox{labelled $n$-shells}}\square[n+1]\fr{}&& K^{(n+1)}. \cocartesien 
}
\] 
Since $(\de\square[n+1])_p = (\square[n+1])_p$ for $p\leq n$, one has
$(K^{(n+1)})_p = (K^{(n)})_p$ for $p\leq n$. And by construction,
$(K^{(n+1)})_{n+1}$ is the set of non-twisted labelled $n$-shell of
$K$.  There is an inclusion map $K^{(n)} \rightarrow K^{(n+1)}$.

\begin{nota} Let $K$ be a $1$-dimensional labelled precubical set with
  $K_0=[p]$ for some $p\geq 2$.  Then let
  \[\COSK(K):= \liminj_{n\geq 1} K^{(n)}.\]
\end{nota}

The important property of the $\COSK$ operator, also used in the proof
of Proposition~\ref{p1eq}, is the following one ensuring that the
higher dimensional automata paradigm is satisfied indeed:

\bp \label{pourquoi_non_degenere} (\cite{ccsprecub} Proposition~3.16)
Let $\square[n]$ be a labelled precubical set with $n\geq 2$. Then one
has the isomorphism of labelled precubical sets
$\COSK(\square[n]_{\leq 1}) \iso \square[n]$.  \ep

Roughly speaking, Proposition~\ref{pourquoi_non_degenere} states that
for all $1\leq p\leq n-1$, there is a bijection between the
non-twisted labelled $p$-shells of $\square[n]$ and the $(p+1)$-cubes
of $\square[n]$.  If the condition non-twisted is removed, i.e. if we
work with a too naive notion of labelled coskeleton construction as in
\cite{exHDA}, then Proposition~\ref{pourquoi_non_degenere} is no
longer true.  Indeed, a naive labelled coskeleton construction adds
too many cubes whereas the higher dimensional automata paradigm states
that one must recover exactly $\square[n]$ from $\square[n]_{\leq 1}$
which corresponds to the concurrent execution of $n$ actions. That was
the problem in K.  Worytkiewicz's coskeleton construction, which was
corrected in the companion paper \cite{ccsprecub}.

\bd Let $K$ and $L$ be two labelled precubical sets. The {\rm
  synchronized tensor product} is by definition
\[K \ot_\sigma L := \liminj_{\square[m]\rightarrow
  K}\liminj_{\square[n]\rightarrow L} \COSK(Z)\] where $Z$ is the
$1$-dimensional precubical set defined by:
\begin{itemize}
\item $Z_0:= \square[m]_0 \p \square[n]_0$ 
\item $Z_1:=(\square[m]_1 \p \square[n]_0) \oplus (\square[m]_0 \p
  \square[n]_1) \oplus \{(x,y)\in \square[m]_1 \p \square[n]_1,
  \overline{\ell(x)} = \ell(y)\}$ with an obvious definition of the
  face maps and the labelling $\tau =\ell(x,y)$ if $\overline{\ell(x)}
  = \ell(y)$.
\end{itemize}
\ed 

\subsection*{Construction of the labelled precubical set of paths}

\bd A labelled precubical set $\ell:K\rightarrow !\Sigma$ {\rm
  decorated by process names} is a labelled precubical set together
  with a set map $d:K_0 \rightarrow \proc_\Sigma$ called the {\rm
  decoration}. \ed

It is recalled in Table~\ref{combsem} the construction of the labelled
precubical set $\square\llbracket P\rrbracket$ of \cite{ccsprecub} by
induction on the syntax of the name.  The labelled precubical set
$\square\llbracket P\rrbracket$ has a unique initial state canonically
decorated by the process name $P$ and its other states will be
decorated as well in an inductive way. Therefore for every process
name $P$, $\square\llbracket P\rrbracket$ is an object of
$\{i\}\ddownarrow \square^{op}\set \ddownarrow !\Sigma$.

\begin{table}
\begin{center}
\begin{tabular}{|l|}
\hline
$\square\llbracket nil\rrbracket:=\square[0]$ \\
\hline
$\square\llbracket \mu.nil\rrbracket:=\mu.nil
\stackrel{(\mu)}\longrightarrow nil$\\
\hline
$\xymatrix{
\square[0]=\{0\} \ar@{->}[r]^-{0\mapsto nil} \ar@{->}[d]^-{0\mapsto P} & \square\llbracket \mu.nil\rrbracket \ar@{->}[d] \\
\square\llbracket P\rrbracket \ar@{->}[r] & \cocartesien {\square\llbracket \mu.P\rrbracket}}$\\
\hline
$\square\llbracket P+Q\rrbracket := \square\llbracket P\rrbracket \oplus \square\llbracket Q\rrbracket$\\ 
with the binary coproduct taken in $\{i\}\ddownarrow \square^{op}\set \ddownarrow !\Sigma$ \\
\hline
$\xymatrix{
  \square\llbracket (\nu a) P\rrbracket \ar@{->}[r] \ar@{->}[d] \cartesien & \square\llbracket P\rrbracket \ar@{->}[d] \\
  !(\Sigma\backslash (\{a,\overline{a}\})) \ar@{->}[r]
  & !\Sigma}$\\
\hline
$\square\llbracket P||Q\rrbracket  := \square\llbracket P\rrbracket \ot_\sigma \square\llbracket Q\rrbracket$ \\ 
\hline
$\square\llbracket \rec(x)P(x)\rrbracket:=\liminj\limits_n \square\llbracket P^n(nil)\rrbracket$\\
\hline
\end{tabular}
\caption{Combinatorial semantics of CCS}
\label{combsem}
\end{center}
\end{table}

\section{Categorical semantics of CCS}
\label{catsem}

The categorical semantics of CCS is obtained from the combinatorial
semantics by applying the geometric realization functor $| - | :
\square^{op}\set \rightarrow \dtop$ introduced in
\cite{ccsprecub}~\footnote{Of course, all theorems proved in the case
  of compactly generated topological spaces in \cite{ccsprecub} are
  still available in the case of $\Delta$-generated topological spaces
  since they only depend on the model structure on $\dtop$.}.  As
already explained in \cite{ccsprecub}, it is necessary to use the
model structure introduced in \cite{model3}, and adapted in
\cite{interpretation-glob} for the framework of $\Delta$-generated
topological spaces. Equivalent geometric realization functors are
defined in \cite{realization}. We will use the construction of
\cite{ccsprecub} in this paper.

Let $Z$ be a topological space. The flow $\glob(Z)$ is defined by
\begin{itemize}
\item $\glob(Z)^0=\{\widehat{0},\widehat{1}\}$, 
\item $\P \glob(Z)= \P_{\widehat{0},\widehat{1}} \glob(Z) = Z$,
\item $s=\widehat{0}$, $t=\widehat{1}$ and a trivial composition law.
\end{itemize}
It is called the \textit{globe} of the space $Z$.

The model structure of \cite{interpretation-glob} is characterized as
follows:
\begin{itemize}
\item The weak equivalences are the \textit{weak S-homotopy
    equivalences}, i.e. the morphisms of flows $f:X\longrightarrow Y$
  such that $f^0:X^0\longrightarrow Y^0$ is a bijection of sets and
  such that $\P f:\P X\longrightarrow \P Y$ is a weak homotopy
  equivalence.
\item The fibrations are the morphisms of flows $f:X\longrightarrow Y$
  such that $\P f:\P X\longrightarrow \P Y$ is a Serre fibration.
\end{itemize}
This model structure is cofibrantly generated. The set of generating
cofibrations is the set $I^{gl}_{+}=I^{gl}\cup
\{R:\{0,1\}\longrightarrow \{0\},C:\varnothing\longrightarrow \{0\}\}$
with
\[\boxed{I^{gl}=\{\glob(\mathbf{S}^{n-1})\subset \glob(\mathbf{D}^{n}), n\geq
0\}}\] where $\mathbf{D}^{n}$ is the $n$-dimensional disk and
$\mathbf{S}^{n-1}$ the $(n-1)$-dimensional sphere. By convention, the
$(-1)$-dimensional sphere is the empty space. The set of generating
trivial cofibrations is
\[\boxed{J^{gl}=\{\glob(\mathbf{D}^{n}\p\{0\})\subset
\glob(\mathbf{D}^{n}\p [0,1]), n\geq 0\}}.\]  

The mapping from $\Obj(\square)$ (the set of objects of $\square$) to
$\Obj(\dtop)$ (the class of flows) defined by $[0] \mapsto \{0\}$ and
$[n] \mapsto \{\widehat{0}<\widehat{1}\}^n$ for $n\geq 1$ induces a
functor from the category $\square$ to the category $\dtop$ by
composition
\[\square \subset \poset \longrightarrow \dtop.\]
\begin{nota}
  A state of the flow $\{\widehat{0}<\widehat{1}\}^n$ is denoted by a
  $n$-uple of elements of $\{\widehat{0},\widehat{1}\}$. The unique
  morphism/execution path from $(x_1,\dots,x_n)$ to $(y_1,\dots,y_n)$
  is denoted by a $n$-uple $(z_1,\dots,z_n)$ with $z_i=x_i$ if
  $x_i=y_i$ and $z_i=*$ if $x_i<y_i$. For example in the flow
  $\{\widehat{0}<\widehat{1}\}^2$ depicted in Figure~\ref{cube2}, one
  has the algebraic relation $(*,*) = (\widehat{0},*)*(*,\widehat{1})
  = (*,\widehat{0}) * (\widehat{1},*)$.
\end{nota}

\bd \cite{ccsprecub} Let $K$ be a precubical set. By definition, the
{\rm geometric realization} of $K$ is the flow
\[\boxed{|K| := \liminj_{\square[n]\rightarrow
    K} (\{\widehat{0}<\widehat{1}\}^n)^{cof}}. \]
\ed

The following proposition is helpful to understand what this geometric
realization functor is. The principle of its proof will be reused in
the paper.

\bp \label{rea-hocolim} Let $K$ be a precubical set. One has a natural
weak S-homotopy equivalence
\[|K| \simeq \holiminj_{\square[n]\rightarrow K}
\{\widehat{0}<\widehat{1}\}^n. \] \ep

\bpf Consider the category of cubes $\square\ddownarrow K$ of $K$. It
is defined by the pullback diagram of small categories
\[
\xymatrix{
\square \ddownarrow K \fr{}\fd{}\cartesien && \square^{op}\set\ddownarrow K
\fd{}\\
&&\\
\square \fr{} && \square^{op}\set.}
\]
In other terms, an object of $\square \ddownarrow K$ is a morphism
$\square[m]\rightarrow K$ and a morphism of $\square\ddownarrow K$ is
a commutative diagram
\[
\xymatrix{
\square[m] \ar@{->}[rd] \ar@{->}[rr] && \square[n] \ar@{->}[ld]\\
& K.&}
\] The category $\square\ddownarrow K$ is a Reedy direct category with
the degree function $d(\square[n]\rightarrow K)=n$.  Since this Reedy
category is direct, the matching category is always empty. So by
\cite{ref_model2} Proposition~15.10.2, the Reedy category has fibrant
constants and the colimit functor $\liminj:\dtop^{\square\ddownarrow
  K} \rightarrow \dtop$ is a left Quillen functor if
$\dtop^{\square\ddownarrow K}$ is equipped with the Reedy model
structure by \cite{ref_model2} Theorem~15.10.8.  Consider the functor
$D:\square\ddownarrow K \rightarrow \dtop$ defined by
$D(\square[n]\rightarrow K) := (\{\widehat{0}<\widehat{1}\}^n)^{cof}$.
One has to check that the diagram $D$ is Reedy cofibrant and the proof
will be complete. By definition of the Reedy model structure, it
suffices to show that for all $n\geq 0$, and with $\alpha =
\square[n]\rightarrow K$, the map $L_{\alpha}D \rightarrow D(\alpha)$
is a cofibration where $L_{\alpha}D$ is the latching object at
$\alpha$. It is easy to see that the latter map is the morphism of
flows $|\de\square[n] \subset \square[n]|$ which is a cofibration of
flows by Theorem~\ref{pasdepb} below.  \epf

\begin{figure}
\[
\xymatrix{
(\widehat{0},\widehat{0}) \fr{(\widehat{0},*)}\fd{(*,\widehat{0})}\ar@{->}[ddrr]^-{(*,*)} && (\widehat{0},\widehat{1})\fd{(*,\widehat{1})}\\
&& \\
(\widehat{1},\widehat{0}) \fr{(\widehat{1},*)}&& (\widehat{1},\widehat{1})}
\]
\caption{The flow $|\square[2]|_{bad} = \{\widehat{0}<\widehat{1}\}^2$
  ($(*,*) = (\widehat{0},*)*(*,\widehat{1}) = (*,\widehat{0}) *
  (\widehat{1},*)$)}
\label{cube2}
\end{figure}

The functor $[n] \mapsto \{\widehat{0}<\widehat{1}\}^n$ from $\square$
to $\dtop$ also induces a \textit{bad} realization functor from
$\square^{op}\set$ to $\dtop$ defined by
\[\boxed{|K|_{bad} := \liminj_{\square[n]\rightarrow K}
  \{\widehat{0}<\widehat{1}\}^n}.\] 
This functor is a bad realization because of the following bad
behaviour:

\bth \label{pbpb} (\cite{ccsprecub} Theorem~7.2) Let $n\geq 3$. The
inclusion of precubical sets $\de\square[n] \subset \square[n]$
induces an isomorphism $|\de\square[n]|_{bad} \iso
|\square[n]|_{bad}$. \eth

On the contrary, the geometric realization functor is well-behaved:

\bth \label{pasdepb} (\cite{ccsprecub} Proposition~7.6 and
\cite{ccsprecub} Theorem~7.8) For any $n\geq 0$, the map of flows
$|\de\square[n] \subset \square[n]|$ is a non-trivial cofibration of
flows.  Moreover, the path space $\P_{\widehat{0}\dots\widehat{0},
  \widehat{1}\dots\widehat{1}}|\de\square[n]|$ is homotopy equivalent
to $\mathbf{S}^{n-2}$ and the path space
$\P_{\widehat{0}\dots\widehat{0},
  \widehat{1}\dots\widehat{1}}|\square[n]|$ is contractible. \eth

Let $K\rightarrow !\Sigma$ be a labelled precubical set. Then the
composition $|K| \rightarrow |!\Sigma| \rightarrow |!\Sigma|_{bad}
\iso ?\Sigma$ gives rise to a labelled flow by \cite{ccsprecub}
Proposition~8.1.

\begin{nota}
  For every process name $P$, let $\llbracket P\rrbracket := |
  \square\llbracket P\rrbracket |$. The flow $\llbracket P\rrbracket$
  is always cofibrant by \cite{ccsprecub} Proposition~7.7.
\end{nota}

Note that a decorated labelled precubical set gives rise to a
decorated labelled flow in the following sense:

\bd A labelled flow $\ell:X \rightarrow ?\Sigma$ {\rm decorated
  by process names} is a labelled flow together with a set map $d:X^0
  \rightarrow \proc_\Sigma$ called the {\rm decoration}. \ed

\section{Relevance of weak S-homotopy for concurrency theory}
\label{relevance}

The translation of the combinatorial semantics of CCS into a
categorical semantics in terms of flows requires the use of
\textit{non-canonical} constructions, more precisely, a non-canonical
choice of a cofibrant replacement functor, and also later
non-canonical choices for homotopy limits and homotopy colimits. The
following theorem is therefore very important:

\bth \label{noncan} For any flow $X$, there exists at most one
precubical set $K$ up to isomorphism such that $|K| \simeq X$. In
other terms, the functor \[K\mapsto \hbox{weak S-homotopy type of
}|K|\] from $\square^{op}\set$ to the homotopy category of flows
$\ho(\dtop)$ reflects isomorphisms. \eth

The precubical set $K$ does not necessarily exist. For example,
$\glob(\mathbf{S}^1)$ is not weakly S-homotopy equivalent to any
geometric realization of any precubical set. Indeed, if there existed
a precubical set $K$ with $|K| \simeq \glob(\mathbf{S}^1)$, then $K$
would have a unique initial state $\widehat{0}$ and a unique final
state $\widehat{1}$, so $K_0 = \{\widehat{0},\widehat{1}\}$. So the
only possibility is a set of $1$-cubes from $\widehat{0}$ to
$\widehat{1}$. Thus the space $\P(|K|)$ would be homotopy equivalent
to a discrete space.

Before proving Theorem~\ref{noncan}, we need to establish several
preliminary results involving among other things the simplicial
structure of the category of flows.

In any flow $X$, if two execution paths $x$ and $y$ are in the same
path-connected component of some $\P_{\alpha,\beta}X$, then there
exists a continuous map $\phi:[0,1] \rightarrow \P_{\alpha,\beta}X$
with $\phi(0) = x$ and $\phi(1) = y$. So for any execution path $z$
such that $x*z$ and $y*z$ exist, the continuous map $\psi$ from
$[0,1]$ to $\P X$ defined by $\psi(t)=\phi(t)*z$ is a continuous path
from $x*z$ to $y*z$. Hence:

\begin{nota} Any flow $X$ induces a flow over the category of sets
  denoted by $\widehat{\pi}_0(X)$ defined by
  $\widehat{\pi}_0(X)^0=X^0$, $\P\widehat{\pi}_0(X)=\pi_0(\P X)$ where
  $\pi_0$ is the path-connected component functor and with the
  composition law induced by the one of $X$. \end{nota}

\begin{nota} Let $\dtop(\set)$ be the category of flows enriched over
  sets, i.e. of small categories without identity maps. \end{nota}

\bp \label{left-adjoint} The functor $\widehat{\pi}_0:\dtop
\rightarrow \dtop(\set)$ is a left adjoint. In particular, it is
colimit-preserving. \ep

\bpf It suffices to prove that the path-connected component functor
$\pi_0:\top \rightarrow \set$ is a left adjoint (let us repeat that we
are working with $\Delta$-generated topological spaces). Here are two
possible arguments:
\begin{enumerate}
\item Every space is homeomorphic to the disjoint sum of its
  path-connected components by \cite{interpretation-glob}
  Proposition~2.8. In fact, a space is even connected if and only if
  it is path-connected. So it is easy to see that the right adjoint is
  the functor from $\set$ to $\top$ taking a set $S$ to the discrete
  space $S$.
\item The functor from $\set$ to $\top$ taking a set $S$ to the
  discrete space $S$ commutes with limits because there is no
  non-discrete totally disconnected $\Delta$-generated spaces, and
  with colimits as in the category of general topological spaces. In
  particular it is accessible. So by \cite{MR95j:18001} Theorem~1.66,
  it has a left adjoint and it is easy to see that the left adjoint is
  the path-connected component functor.
\end{enumerate}
\epf 

\bp \label{P-accessible} (Compare with \cite{interpretation-glob}
Proposition~4.9) The path space functor $\P:\dtop \rightarrow \top$ is
a right adjoint. In particular, it is accessible.  \ep

In fact, the functor $\P:\dtop \rightarrow \top$ is of course finitely
accessible.

\bpf Let $Z$ be a topological space. By \cite{interpretation-glob}
Proposition~2.8, $Z$ is homeomorphic to the disjoint union of its
path-connected components. Let us write this situation by \[Z\iso
\bigsqcup_{Z_i\in\pi_0(Z)} Z_i.\] Then one has for any flow $X$ \beas
\top(Z,\P X) &\iso& \prod_{Z_i\in\pi_0(Z)} \top(Z_i,\P X) \\
&\iso& \prod_{Z_i\in\pi_0(Z)} \dtop(\glob(Z_i),X) \\
&\iso& \dtop(\bigsqcup_{Z_i\in\pi_0(Z)}\glob(Z_i),X).  \eeas So the
path space functor $\P:\dtop \rightarrow \top$ is accessible by
\cite{MR95j:18001} Theorem~1.66.  \epf

\bp \label{inc-path} Let $i:A\rightarrow X$ be a cofibration of flows
between cofibrant flows. Then the continuous map $\P i:\P A
\rightarrow \P X$ is a cofibration between cofibrant spaces.  \ep

Note that Proposition~\ref{inc-path} remains true if we only suppose
that the space $\P A$ is cofibrant. Proposition~\ref{inc-path} is a
generalization of \cite{interpretation-glob} Proposition~7.5.

\bpf Let us suppose first that there is a pushout diagram of flows
\[
\xymatrix{
\glob(\mathbf{S}^{n-1}) \fr{} \fd{} && A \fd{i} \\
&&\\
\glob(\mathbf{D}^{n}) \fr{} && X.\cocartesien}
\] 
By \cite{model3} Proposition~15.1, the continuous map $\P i:\P A
\rightarrow \P X$ is a transfinite composition of pushouts of maps of
the form $\id_{X_1}\p\dots\p i_n\p \dots \p \id_{X_p}$ where the
spaces $X_i$ are spaces of the form $\P_{\alpha,\beta}A$ and where
$i_n:\mathbf{S}^{n-1}\subset \mathbf{D}^n$ is the inclusion with
$n\geq 0$. Any space of the form $\P_{\alpha,\beta}A$ is cofibrant by
\cite{interpretation-glob} Proposition~7.5 since $A$ is cofibrant.  So
the map $\P i:\P A \rightarrow \P X$ is a cofibration because the
model category $(\top,\p)$ is monoidal.

Let us treat now the general case. The cofibration $i$ is a retract of
a map $j:A\rightarrow Y$ of $\cell(I^{gl}_{+})$ by a map which fixes
$A$ by \cite{MR99h:55031} Corollary~2.1.15. So the continuous map $\P
i : \P A \rightarrow \P X$ is a retract of the continuous map $\P j :
\P A \rightarrow \P Y$.  The map of flows $j : A \rightarrow Y$ is the
composition of a transfinite sequence $Z : \lambda \rightarrow \dtop$
for some ordinal $\lambda$ with $Z_0 = A$. By
Proposition~\ref{P-accessible}, one has the homeomorphism $\liminj \P
Z_\alpha \iso \P Y$. The first part of this proof implies that $\P j :
\P A \rightarrow \P Y$ is then a cofibration of spaces, and therefore
that $\P i:\P A \rightarrow \P X$ is a cofibration as well.  \epf

\begin{nota} The associative monoid without unit $(\N^{*},+)$ of
  strictly positive integers together with the addition can be viewed
  as a flow with one object, the discrete path space $\N^{*}$ and the
  composition law $+$. \end{nota}

\begin{nota} Let $K$ be a precubical set. Let $K_{\leq n}$ be the
  precubical set obtained from $K$ by keeping the $p$-dimensional
  cubes of $K$ only for  $p\leq n$. In particular, $K_{\leq 0}=K_0$.
\end{nota}

\bp \label{length} Let $K$ be a precubical set. There exists a unique
morphism of flows $L_K:\widehat{\pi}_0(|K|) \rightarrow \N^{*}$,
natural with respect to $K$, such that for any $x \in K_1$, for any $z
\in \P(|\square[1]|)$, one has $L_K(|x|(z)) = 1$.  \ep

\bpf We construct $L_K:\widehat{\pi}_0(|K_{\leq n}|) \rightarrow
\N^{*}$ for any precubical set $K$ by induction on $n\geq 0$. There is
nothing to do for $n=0$. The passage from $|K_{\leq n}|$ to $|K_{\leq
  n+1}|$ is done as usual by the following pushout diagram of flows:
\[
\xymatrix{
\bigsqcup_{x:\de\square[n+1]\rightarrow K} |\de\square[n+1]|
\fd{}\fr{} && |K_{\leq n}|\fd{}\\
&&\\
\bigsqcup_{x:\de\square[n+1]\rightarrow K} |\square[n+1]|
\fr{}&& \cocartesien {|K_{\leq n+1}|}}
\] 
where the sum is over all $n$-shells $x:\de\square[n+1]\subset
\square[n+1]\rightarrow K$. Let $n\geq 0$. By induction hypothesis,
the flow $\widehat{\pi}_0(|\de\square[n+1]|)$ and
$\widehat{\pi}_0(|K_{\leq n}|)$ are defined. We know that the map of
flows $|\de\square[n+1]| \rightarrow |\square[n+1]|$ is a cofibration
by Theorem~\ref{pasdepb}.  In fact, this map of flows induces the
identity maps $\P_{\alpha,\beta} (|\de\square[n+1]|) =
\P_{\alpha,\beta} (|\square[n+1]|)$ for $(\alpha,\beta) \neq
(\widehat{0}\dots\widehat{0}, \widehat{1}\dots\widehat{1})$ and a
non-trivial cofibration~\footnote{Let us recall that the space
  $\P_{\widehat{0}\dots\widehat{0},\widehat{1}\dots\widehat{1}}(|\square[n+1]|)$
  is contractible and that by Theorem~\ref{pasdepb}, there is a
  homotopy equivalence
  $\P_{\widehat{0}\dots\widehat{0},\widehat{1}\dots\widehat{1}}(|\de\square[n+1]|)
  \simeq \mathbf{S}^{n-1}$.} between cofibrant spaces
$\P_{\widehat{0}\dots\widehat{0},\widehat{1}\dots\widehat{1}}(|\de\square[n+1]|)
\rightarrow
\P_{\widehat{0}\dots\widehat{0},\widehat{1}\dots\widehat{1}}(|\square[n+1]|)$
by Proposition~\ref{inc-path}. Then let $L_{\square[n+1]}(x)=n+1$ for
any $x \in
\P_{\widehat{0}\dots\widehat{0},\widehat{1}\dots\widehat{1}}(|\square[n+1]|)$.
One obtains the commutative square of $\dtop(\set)$:
\[
\xymatrix{ \bigsqcup_{x:\de\square[n+1]\rightarrow K}
  \widehat{\pi}_0(|\de\square[n+1]|)
  \fd{}\fr{} && \widehat{\pi}_0(|K_{\leq n}|)\fd{}\\
  &&\\
  \bigsqcup_{x:\de\square[n+1]\rightarrow K} \widehat{\pi}_0(|\square[n+1]|)
  \fr{}&& {\N^{*}}.}
\] 
By Proposition~\ref{left-adjoint} and by the universal property of the
pushout, one obtains the natural map $\widehat{\pi}_0(|K_{\leq n+1}|)
\rightarrow \N^{*}$. Since the functor $K\mapsto \widehat{\pi}_0(|K|)$ is a left
adjoint, one obtains a natural map $\widehat{\pi}_0(|K|) \iso \liminj
\widehat{\pi}_0(|K_{\leq n}|) \rightarrow \N^{*}$. \epf

\bd The integer $L_K(x)$ for $x\in \P(|K|)$ is called the {\rm length}
of $x$. \ed

Proposition~\ref{length} means that the length of $x\in \P(|K|)$
satisfies the following (intuitive) algebraic rules: 
\begin{itemize}
\item $L_K(x*y) = L_K(x) + L_K(y)$ if $x$ and $y$ are composable
\item $L_K(x) = L_K(y)$ if $x$ and $y$ are in the same path-connected
  component of the space $\P(|K|)$
\item $L_K(x)=1$ if $x$ corresponds to an edge, i.e. a $1$-cube, of
  the precubical set $K$
\item the naturality of the morphism of flows $L_K :
  \widehat{\pi}_0(|K|) \rightarrow \N^{*}$ means that length is
  preserved by a map of precubical sets.
\end{itemize}

The model category $\dtop$ is simplicial by \cite{realization}
Section~3 and \cite{interpretation-glob} Appendix~B. Let $\map(X,Y)$
be the function complex from $X$ to $Y$. It is equal to the simplicial
nerve of the space $\tdtop(X,Y)$ of morphisms of flows from $X$ to $Y$
equipped with the Kelleyfication of the relative topology.

\bp \label{inj} Let $K$ be a precubical set.  Let $n\geq 0$. The
natural set map $K_n \rightarrow \dtop(|\square[n]|,|K|)$ defined by
taking $x\in K_n$ to $|x|:|\square[n]| \rightarrow |K|$ is one-to-one.
\ep

\bpf One has $|K| = \liminj_{\square[n]\rightarrow K} |\square[n]|$ by
definition. If $x$ and $y$ are two different $n$-cubes of $K$, then
they correspond to two different copies of $|\square[n]|$ in the
colimit calculating $K$.  Let \[\gamma\in
\P_{\widehat{0}\dots\widehat{0},\widehat{1}\dots\widehat{1}}|\square[n]|
\backslash
\P_{\widehat{0}\dots\widehat{0},\widehat{1}\dots\widehat{1}}|\de\square[n]|.\]
Then $|x|(\gamma) \neq |y|(\gamma)$. Therefore $|x| \neq |y|$.  \epf

\begin{nota} Let $K$ be a precubical set. The precubical set $\widehat{K}$ is
defined by \[\boxed{\widehat{K} = \pi_0 \map(|\square[*]|,|K|) =
  \pi_0 \tdtop(|\square[*]|,|K|)}.\]
\end{nota}

Since $|\square[n]|$ is cofibrant and since all flows are fibrant, the
function complex $\map(|\square[n]|,|K|)$ is weakly equivalent to the
homotopy function complex from $|\square[n]|$ to $|K|$. Thus
$\widehat{K}_n = \ho(\dtop)(|\square[n]|,|K|)$ for all $n\geq 0$ where
$\ho(\dtop)$ is the homotopy category of $\dtop$.

The natural map of precubical sets $K \rightarrow
\dtop(|\square[*]|,|K|)$ induces a natural map of precubical sets
$K\rightarrow \widehat{K}$.

\bp \label{inj2} Let $K$ be a precubical set.  Let $n\geq 0$. The
continuous map $j_n:\tdtop(|\square[n]|,|K_{\leq n}|) \rightarrow
\tdtop(|\square[n]|,|K|)$ induced by the inclusion of precubical sets
$K_{\leq n} \subset K$ is an inclusion of $\Delta$-generated spaces in
the sense that one has a homeomorphism \[\tdtop(|\square[n]|,|K_{\leq
  n}|) \iso j_n(\tdtop(|\square[n]|,|K_{\leq n}|))\] with the
right-hand topological space equipped with the Kelleyfication of the
relative topology. \ep

\bpf[Sketch of proof] The map $j_n:\tdtop(|\square[n]|,|K_{\leq n}|)
\rightarrow \tdtop(|\square[n]|,|K|)$ is clearly one-to-one.  It
suffices to prove that for any continuous map $\phi : Z \rightarrow
\tdtop(|\square[n]|,|K|)$ such that $\phi(Z)\subset
j_n(\tdtop(|\square[n]|,|K_{\leq n}|))$, the unique set map
$Z\rightarrow \tdtop(|\square[n]|,|K_{\leq n}|)$ induced by $\phi$ is
continuous.

By Theorem~\ref{pasdepb}, the map $|K_{\leq n}| \rightarrow |K|$ is a
cofibration of flows. One has $|K_{\leq n}|^0 = |K|^0 = K_0$ and the
continuous map $\P(|K_{\leq n}|) \rightarrow \P(|K|)$ is a cofibration
of spaces by Proposition~\ref{inc-path}. So the latter continuous map
is a closed $T_1$-inclusion of general topological spaces by
\cite{MR99h:55031} Lemma~2.4.5, and also an inclusion of
$\Delta$-generated spaces.

By \cite{interpretation-glob} Appendix~B, the category of flows
enriched over $\Delta$-generated topological spaces is tensored and
cotensored over the $\Delta$-generated spaces in the sense of
\cite{strom2}. So the continuous map $\phi : Z \rightarrow
\tdtop(|\square[n]|,|K|)$ corresponds by adjunction to a morphism of
flows $|\square[n]| \ot Z \rightarrow |K|$.  By hypothesis, the map
$\P \phi$ factors uniquely as a set map as a composite
\[\P(|\square[n]| \ot Z) \rightarrow \P(|K_{\leq n}|) \rightarrow
\P(|K|).\] Since the right-hand map is a closed $T_1$-inclusion of
general topological spaces, the left-hand map $\P(|\square[n]| \ot Z)
\rightarrow \P(|K_{\leq n}|)$ is continuous.  Hence the factorization
$|\square[n]| \ot Z \rightarrow |K_{\leq n}| \rightarrow |K|$. By
adjunction, one obtains the continuous map $Z\rightarrow
\tdtop(|\square[n]|,|K_{\leq n}|)$.  \epf

\bp \label{reflete} The functor $K\mapsto \widehat{K}$ reflects
isomorphisms, i.e a map of precubical sets $f:K \rightarrow L$ is an
isomorphism if and only if the map of precubical sets
$\widehat{f}:\widehat{K} \rightarrow \widehat{L}$ is an isomorphism.
\ep

\bpf It turns out that the natural map of precubical sets
$K\rightarrow \widehat{K}$ is a monomorphism. Indeed, take two
elements $x$ and $y$ of $K_n$ such that $|x|$ and $|y|$ are in the
same path-connected component of $\tdtop(|\square[n]|,|K|)$. By
definition, there exists a continuous map $\phi : [0,1] \rightarrow
\tdtop(|\square[n]|,|K|)$ such that $\phi(0) = |x|$ and $\phi(1) =
|y|$. For any $z\in \P (|\square[n]|)$, one has the inequality
$L_K(\phi(t)(z)) \leq n$ for all $t\in [0,1]$ because
$L_{\square[n]}(z) \leq n$ and because maps of precubical sets
preserve length.  But any execution path of $\P(|K|)\backslash
\P(|K_{\leq n}|)$ is of length strictly greater than $n$. So the map
$\phi$ factors uniquely as a composite $[0,1] \rightarrow
\tdtop(|\square[n]|,|K_{\leq n}|) \rightarrow
\tdtop(|\square[n]|,|K|)$ by Proposition~\ref{inj2}. Since a
non-trivial homotopy $\phi$ would necessarily use higher dimensional
cubes of $K\backslash K_{\leq n}$, the homotopy $\phi$ is trivial.
Therefore $|x| = |y|$, and by Proposition~\ref{inj} one obtains $x =
y$.

So the precubical set $K$ is naturally isomorphic to a precubical
subset of $\widehat{K}$. Take a map $f:K\rightarrow L$. Then, by
naturality, there is a commutative square of precubical sets
\[
\xymatrix{
K\fd{f} \fr{} && \widehat{K} \fd{\widehat{f}}\\
&&\\
L \fr{}&& \widehat{L}.}
\] 
If $f$ is not an isomorphism, then two situations may happen:
\begin{itemize}
\item There exist $n\geq 0$ and two distinct $n$-cubes $x$ and $y$ of
  $K$, and therefore of $\widehat{K}$, with $f(x) = f(y)$. Then
  $\widehat{f}(x) = \widehat{f}(y)$ and therefore $\widehat{f}$ is not
  an isomorphism.
\item There exist $n\geq 0$ and a $n$-cube $x$ of $L$ which does not
  belong to the image of $f$. Since the map $\widehat{f}$ factors as a
  composite $\widehat{K} \rightarrow \widehat{f(K)} \rightarrow
  \widehat{L}$, the $n$-cube $x$ does not have any antecedent by
  $\widehat{f}$. So $\widehat{f}$ is not an isomorphism.
\end{itemize}
\epf 

\bpf[Proof of Theorem~\ref{noncan}] Let $K$ and $L$ be two precubical
sets with $|K|\simeq |L|$. For all $n\geq 0$, the functor
$\map(|\square[n]|,-):\dtop \rightarrow \sis$ preserves weak
equivalences between fibrant objects by \cite{ref_model2}
Corollary~9.3.3 since this functor is a right Quillen functor. So
there is an isomorphism $\widehat{K}\iso \widehat{L}$ since both $|K|$
and $|L|$ are fibrant~\footnote{All flows are actually fibrant.}. And
by Proposition~\ref{reflete}, one obtains an isomorphism $K\iso L$.
\epf

In conclusion, we can safely work up to weak S-homotopy without losing
any kind of computer-scientific information already present in the
structure of the precubical set.

\section{Effect of the geometric realization functor when it is a left
adjoint}
\label{effect}

One has the isomorphism $\llbracket P+Q\rrbracket \iso \llbracket P\rrbracket \oplus \llbracket Q\rrbracket$ of
$\{i\}\ddownarrow \dtop \ddownarrow !\Sigma$ since the geometric
realization functor is a left adjoint.

\bp \label{L1} One has the pushout diagram of labelled flows
\[
\xymatrix{
  \{0\} \fr{0\mapsto nil} \fd{0\mapsto i} && \llbracket \mu.nil\rrbracket \fd{} \\
  &&\\
  \llbracket P\rrbracket \fr{} && \cocartesien {\llbracket \mu.P\rrbracket}}\] and this diagram is also
a homotopy pushout diagram.  \ep

\bpf The diagram above is a pushout diagram since the geometric
realization functor is a left adjoint. This diagram is also a homotopy
pushout diagram by \cite{MR99h:55031} Lemma~5.2.6 since the three
flows $\{0\}$, $\llbracket \mu.nil\rrbracket$ and $\llbracket
P\rrbracket$ are cofibrant and since the map $\{0\} \rightarrow
\llbracket P\rrbracket$ is a cofibration.  \epf

\bp \label{L2} Let $P(x)$ be a process name with one free guarded
variable $x$.  Then one has the isomorphism
\[ \llbracket \rec(x)P(x)\rrbracket \iso \liminj_n \llbracket P^n(nil)\rrbracket\] 
and the colimit is also a homotopy colimit. 
\ep

\bpf The isomorphism comes again from the fact that the geometric
realization functor is a left adjoint.  The tower of flows $n\mapsto
\llbracket P^n(nil)\rrbracket$ is a tower of cofibrant flows and each
map $\llbracket P^n(nil)\rrbracket \rightarrow \llbracket
P^{n+1}(nil)\rrbracket$ is a cofibration by Theorem~\ref{pasdepb}. So
the colimit is also a homotopy colimit by \cite{ref_model2}
Proposition~15.10.12.  \epf

\bp \label{unpullback} Let $K \rightarrow !\Sigma$ be a labelled
precubical set.  Let $\Sigma' \subset \Sigma$.  Consider the pullback
diagram of precubical sets
\[
\xymatrix{
L \cartesien \fr{} \fd{} && K \fd{} \\
&& \\
!\Sigma' \fr{} && !\Sigma.}
\] 
Then the commutative diagram of flows 
\[
\xymatrix{
|L| \fr{} \fd{} && |K| \fd{} \\
&& \\
?\Sigma' \fr{} && ?\Sigma}
\] 
obtained by taking the realization of the first diagram and by
composing with the commutative square
\[
\xymatrix{
|!\Sigma'| \fr{} \fd{} && |!\Sigma| \fd{} \\
&& \\
?\Sigma' \fr{} && ?\Sigma}
\] 
is a pullback and a homotopy pullback diagram of flows.  \ep

\bpf It is well-known that every precubical set $K$ is a
$\{\de\square[n] \subset \square[n],n\geq 0\}$-cell complex since the
passage from $K_{\leq n-1}$ to $K_{\leq n}$ for $n\geq 1$ is done by
the following pushout diagram:
\[
\xymatrix{
\bigsqcup_{x\in K_n} \de\square[n] \fr{} \fd{} && K_{\leq n-1} \fd{} \\
&& \\
\bigsqcup_{x\in K_n} \square[n] \fr{} && K_{\leq n}\cocartesien}
\] 
where the map $\de\square[n] \rightarrow K_{\leq n-1}$ indexed by
$x\in K_n$ is induced by the $(n-1)$-shell $\de\square[n] \subset
\square[n] \stackrel{x}\rightarrow K$. One also has the pullback
diagram of sets
\[
\xymatrix{
  L_n \iso \square^{op}\set(\square[n],L) \cartesien \fr{} \fd{} && K_n \iso \square^{op}\set(\square[n],K) \fd{} \\
  && \\
  (!\Sigma')_n \iso \square^{op}\set(\square[n],!\Sigma') \fr{} &&
  (!\Sigma)_n \iso \square^{op}\set(\square[n],!\Sigma)}
\] 
by the Yoneda lemma and because pullbacks are calculated pointwise in
the category of precubical sets. So the precubical set $L$ is the
$\{\de\square[n] \subset \square[n],n\geq 0\}$-cell subcomplex
obtained by keeping the cells $\de\square[n] \subset \square[n]$
induced by the $n$-dimensional cubes $\square[n] \rightarrow K$ such
that the composite $\square[n] \rightarrow K \rightarrow !\Sigma$
factors as a composite $\square[n] \rightarrow !\Sigma' \rightarrow
!\Sigma$. Thus, the map $L \rightarrow K$ is a relative
$\{\de\square[n] \subset \square[n],n\geq 0\}$-cell complex. One has
the bijection $( !\Sigma')_0 \iso ( !\Sigma)_0$.  Therefore $L_0\iso
K_0$ and the map $L \rightarrow K$ is a relative $\{\de\square[n]
\subset \square[n],n\geq 1\}$-cell complex.  Since the realization
functor $K \mapsto |K|$ is a left adjoint, the map $|L| \rightarrow
|K|$ is then a relative $\{|\de\square[n]| \subset |\square[n]|,n\geq
1\}$-cell complex. By Theorem~\ref{pasdepb}, we deduce that the map
$|L| \rightarrow |K|$ is a cofibration of flows with $|L|^0=|K|^0$. By
Proposition~\ref{inc-path}, the continuous map $\P(|L|) \rightarrow
\P(|K|)$ is a [closed $T_1$-]inclusion of general topological spaces
in the sense that for any continuous map $f:Z\rightarrow \P(|K|)$ such
that $f(Z)$ is in the image of $\P(|L|)$, there exists a unique
continuous map $\overline{f}:Z \rightarrow \P(|L|)$ such that the
composition $Z \rightarrow \P(|L|) \rightarrow \P(|K|)$ is equal to
$f$. Consider a commutative diagram of flows
\[
\xymatrix{
W \ar@{-->}[rd]^-{k}\ar@/^10pt/[rrrd]^-{u}\ar@/_10pt/[dddr]_-{v}& && \\
&|L| \fr{} \fd{} && |K| \fd{\ell}  \\
&&& \\
&?\Sigma' \fr{} && ?\Sigma
}
\] 
Let $\gamma\in\P W$. By definition, one has
\[|K| =\liminj\limits_{\square[n] \rightarrow K}
(\{\widehat{0}<\widehat{1}\}^n)^{cof}\] with one copy of
$(\{\widehat{0}<\widehat{1}\}^n)^{cof}$ corresponding to one element
$x \in K_n$.  Thus, $u(\gamma) = \gamma_1 *\dots * \gamma_r$ with
$\gamma_i\in \P (\{\widehat{0}<\widehat{1}\}^{n_i})^{cof}$
corresponding to a $n_i$-dimensional cube $x_i$ of $K$. And
$\ell(\gamma_1 *\dots * \gamma_r) = a_1*\dots*a_s$ with $a_i\in
\Sigma'$ for all $i=1,\dots,s$ (note $r$ is not necessarily equal to
$s$). By construction of $L$, the $n_i$-dimensional cube $x_i$ of $K$
then belongs to $L$. By definition, one has
\[|L| =\liminj\limits_{\square[n] \rightarrow K}
(\{\widehat{0}<\widehat{1}\}^n)^{cof}\] with one copy of
$(\{\widehat{0}<\widehat{1}\}^n)^{cof}$ corresponding to one element
$x \in L_n$. So $u(\gamma)$ belongs to the image of the inclusion of
spaces $\P (|L|) \rightarrow \P (|K|)$. Hence the existence and the
uniqueness of $k$. So the commutative square
\[
\xymatrix{
|L| \fr{} \fd{} && |K| \fd{} \\
&& \\
?\Sigma' \fr{} && ?\Sigma}
\] 
is a pullback diagram of flows. A map of flows $f:X \rightarrow Y$ is
a fibration if and only if the continuous map $\P f:\P X \rightarrow
\P Y$ is a Serre fibration. Therefore all objects of $\dtop$ are
fibrant. And the map $|K| \rightarrow ?\Sigma$ is a fibration of flows
since the path space $\P(?\Sigma)$ is discrete. Thus, the pullback
diagram
\[
\xymatrix{
|L| \cartesien \fr{} \fd{} && |K| \fd{} \\
&& \\
?\Sigma' \fr{} && ?\Sigma}
\] 
is also a homotopy pullback diagram of flows by e.g.
\cite{MR99h:55031} Lemma~5.2.6.  \epf

\begin{cor} \label{L3} Let $P$ be a process name. Then the commutative
  diagram
\[
\xymatrix{
  \llbracket (\nu a)P\rrbracket \fr{} \fd{} && \llbracket P\rrbracket \fd{} \\
  && \\
  ?(\Sigma\backslash (\{a,\overline{a}\})) \fr{} && ?\Sigma}
\] 
is both a pullback diagram and a homotopy pullback diagram of flows.
\end{cor}

The following proposition is crucial to get rid of the coskeleton
construction in the interpretation of the parallel composition with
synchronization.

\bp \label{p1eq} Let $\square[m]$ be a labelled $m$-cube with $m\geq
0$. Let $\square[n]$ be a labelled $n$-cube with $n\geq 0$. Then the
map $|\square[m]\ot_\sigma \square[n]| \rightarrow
|\square[m]\ot_\sigma \square[n]|_{bad}$ is a trivial fibration of
flows.  \ep

\bpf By Theorem~\ref{pbpb} saying that $|\de\square[n]|_{bad}\iso
|\square[n]|_{bad}$ for $n\geq 3$, and since the bad geometric
realization is a left adjoint, one has the pushout diagram of flows:
\[
\xymatrix{
\bigsqcup\limits_{\hbox{labelled $1$-shells}} |\de\square[2]|_{bad}\fd{}  \fr{} && |(\square[m]\ot_\sigma \square[n])_{\leq 1}|_{bad}\fd{} \\
&& \\
\bigsqcup\limits_{\hbox{labelled $1$-shells}} |\square[2]|_{bad} \fr{} && \cocartesien {|\square[m]\ot_\sigma \square[n]|_{bad}}}
\] 
The path space $\P(|(\square[m]\ot_\sigma \square[n])_{\leq
  1}|_{bad})$ contains the free compositions of (composable) $1$-cubes
of $\square[m]\ot_\sigma \square[n]$. The effect of the map
$\P(|(\square[m]\ot_\sigma \square[n])_{\leq 1}|_{bad}) \rightarrow
\P(|\square[m]\ot_\sigma \square[n]|_{bad})$ is to add algebraic
relations $v*w=x*y$ whenever $\ell(v)=\ell(y)$, $\ell(w)=\ell(x)$ and
$\ell(v) * \ell(w) = \ell(w) * \ell(v)$.

The map $|\square[m]\ot_\sigma \square[n]| \rightarrow
|\square[m]\ot_\sigma \square[n]|_{bad}$ induces a bijection
$|\square[m]\ot_\sigma \square[n]|^0 \iso |\square[m]\ot_\sigma
\square[n]|^0_{bad}$. The continuous map $\P(|\square[m]\ot_\sigma
\square[n]|) \rightarrow \P(|\square[m]\ot_\sigma \square[n]|_{bad})$
is a Serre fibration since the space $\P(|\square[m]\ot_\sigma
\square[n]|_{bad})$ is discrete.  Therefore, it remains to prove that
the fibre of the fibration $\P (|\square[m]\ot_\sigma \square[n]|)
\rightarrow \P(|\square[m]\ot_\sigma \square[n]|_{bad})$ over
$x_1*\dots*x_r \in \P(|\square[m]\ot_\sigma \square[n]|_{bad})$ where
$x_1, \dots, x_r \in (\square[m]\ot_\sigma \square[n])_1$ is
contractible. Since the labels of $x_1,\dots,x_r$ commute with one
another~\footnote{\label{moregeneral}For more general synchronization
  algebras, it is not true that all the labels necessarily commute
  with one another.  One has first to set $x_1*\dots *x_r = y_1 *
  \dots * y_s$ where the labels contained in each $y_i$ commute with
  one another and one has then to say that the fibre over $x_1*\dots
  *x_r$ is the product of the contractible fibres over the $y_i$.},
this fibre is equal to the path space
$\P_{\widehat{0}\dots\widehat{0},\widehat{1}\dots\widehat{1}}(|\COSK(\square[r]_{\leq
  1})|)$ of execution paths from the initial state to the final state
of the $r$-cube filled out by the $\COSK$ operator.  So the fibre is
contractible by Proposition~\ref{pourquoi_non_degenere}.  \epf

\bth \label{parho} Let $P$ and $Q$ be two process names of
$\proc_\Sigma$. Then the flow associated with the process $P||Q$ is
weakly S-homotopy equivalent to the flow
\[\holiminj_{\square[m]\rightarrow
  \square\llbracket P\rrbracket}\holiminj_{\square[n]\rightarrow \square\llbracket Q\rrbracket}
|(\square[m]\ot_\sigma \square[n])_{\leq 2}|_{bad}.\] \eth 

Note that by the Fubini theorem for homotopy colimits (e.g.,
\cite{monographie_hocolim} Theorem~24.9) the order of homotopy
colimits is not important.

\bpf[Sketch of proof] By Proposition~\ref{p1eq} and
Theorem~\ref{pbpb}, one has a weak S-homotopy equivalence
\[\holiminj_{\square[m]\rightarrow
  \square\llbracket P\rrbracket}\holiminj_{\square[n]\rightarrow \square\llbracket Q\rrbracket}
|\square[m]\ot_\sigma \square[n]| \stackrel{\simeq}\longrightarrow
\holiminj_{\square[m]\rightarrow
  \square\llbracket P\rrbracket}\holiminj_{\square[n]\rightarrow \square\llbracket Q\rrbracket}
|(\square[m]\ot_\sigma \square[n])_{\leq 2}|_{bad}.\] For similar
reasons to the proof of Proposition~\ref{rea-hocolim}, the double
colimit
\[\liminj_{\square[m]\rightarrow
  \square\llbracket P\rrbracket}\liminj_{\square[n]\rightarrow \square\llbracket Q\rrbracket}
  |\square[m]\ot_\sigma \square[n]|\] is a homotopy colimit because
  the diagram is Reedy cofibrant over a fibrant constant Reedy
  category.  So the canonical map
\[\holiminj_{\square[m]\rightarrow
  \square\llbracket P\rrbracket}\holiminj_{\square[n]\rightarrow \square\llbracket Q\rrbracket}
|\square[m]\ot_\sigma \square[n]|
\stackrel{\simeq}\longrightarrow\liminj_{\square[m]\rightarrow
  \square\llbracket P\rrbracket}\liminj_{\square[n]\rightarrow \square\llbracket Q\rrbracket}
  |\square[m]\ot_\sigma \square[n]|\] is a weak S-homotopy
  equivalence.  Since the geometric realization functor is a left
  adjoint, the right-hand double colimit is isomorphic to
\[\left|\liminj_{\square[m]\rightarrow
    \square\llbracket P\rrbracket}\liminj_{\square[n]\rightarrow \square\llbracket Q\rrbracket}
  \square[m]\ot_\sigma \square[n]\right|,\] hence the result by
\cite{ccsprecub} Proposition~4.6 saying that the operator $\ot_\sigma$
preserves colimits.  \epf

The flow $|(\square[m]\ot_\sigma \square[n])_{\leq 2}|_{bad}$ is
obtained from the flow $|(\square[m]\ot_\sigma \square[n])_{\leq
  1}|_{bad}$ by adding an algebraic rule $x*y=z*t$ for each $4$-uple
$(x,y,z,t)$ such that $\ell(x) = \ell(t)$, $\ell(y) = \ell(z)$ and
$\ell(x)*\ell(y) = \ell(z)* \ell(t)$. So the coskeletal approach has
totally disappeared in the statement of Theorem~\ref{parho}.

\begin{cor}
  Let $P$ and $Q$ be two process names of $\proc_\Sigma$. Then the
  flow associated with the process $P||Q$ is weakly S-homotopy
  equivalent to the flow
\[\liminj_{\square[m]\rightarrow
  \square\llbracket P\rrbracket}\liminj_{\square[n]\rightarrow \square\llbracket Q\rrbracket}
(|(\square[m]\ot_\sigma \square[n])_{\leq 2}|_{bad})^{cof}.\]
\end{cor}

\bpf In the model category of flows, the class of cofibrations which
are monomorphisms is closed under pushout and transfinite composition.
Therefore the cofibrant replacement of a monomorphism is a
cofibration, and even an inclusion of subcomplexes (\cite{ref_model2}
Definition~10.6.7) because the cofibrant replacement functor is
obtained by the small object argument, starting from the identity of
the initial object, i.e. the empty flow. So the diagram calculating
\[\liminj_{\square[m]\rightarrow
  \square\llbracket P\rrbracket}\liminj_{\square[n]\rightarrow \square\llbracket Q\rrbracket}
(|(\square[m]\ot_\sigma \square[n])_{\leq 2}|_{bad})^{cof}\] is Reedy
cofibrant. Thus the double colimit above has the correct weak
S-homotopy type.  \epf

\section{Towards a pure homotopical semantics of CCS}
\label{bonus}

\begin{table}
\begin{center}
\begin{tabular}{|l|}
\hline
$[P]:=|\llbracket P\rrbracket|$ for $P\in\proc_\Sigma$\\
\hline
$\xymatrix{
[\{0\}] \ar@{->}[r]^-{0\mapsto nil} \ar@{->}[d]^-{0\mapsto P} & [\mu.nil] \ar@{->}[d] \\
[P] \ar@{->}[r] & \wcocartesien {[\mu.P]}}$\\
\hline
$[P+Q] := [P] \oplus [Q]$\\ 
with the binary coproduct taken in $\ho(\{i\}\ddownarrow \dtop \ddownarrow ?\Sigma)$ \\
\hline
$\xymatrix{
  [(\nu a) P] \ar@{->}[r] \ar@{->}[d] \wcartesien & [P] \ar@{->}[d] \\
  [?(\Sigma\backslash (\{a,\overline{a}\}))] \ar@{->}[r]
  & [?\Sigma]}$\\
\hline
$[\rec(x)P(x)]:=\wliminj\limits_n [P^n(nil)]$\\
\hline
\end{tabular}
\caption{Pure homotopical semantics of a restriction of CCS, w meaning
  Heller's privileged weak (co)limits of $\ho(\dtop)$}
\label{wsem}
\end{center}
\end{table}

Let us restrict our attention to CCS without parallel composition with
synchronization.  So the new syntax of the language \textit{for this
  section only} is:
\[
P::=P\in\proc_\Sigma \ |\ a.P \ |\ (\nu a)P \ |\ P + P \ |\ \rec(x)P(x).\] Denote
by $\ho(\dtop)$ the homotopy category of flows, i.e. the categorical
localization of the flows by the weak S-homotopy equivalences.  We
want to explain in this section how it is possible to construct a
semantics of this restriction of CCS in terms of elements of
$\ho(\dtop)$. 

The following theorem is about realization of homotopy commutative
diagrams in the particular case of a diagram over a Reedy category. It
gives a sufficient condition for a homotopy commutative diagram to be
coherently homotopy commutative.

\bth (Cisinski) \label{weak} (\cite{CiSi} for the finite case and
\cite{cofibrationcat} Theorem~8.8.5 for the generalization) Let
$\mathcal{M}$ be a model category. Let $\mathcal{B}$ be a small Reedy
category which is free, i.e. freely generated by a graph.  Moreover,
let us suppose that $\mathcal{B}$ is either direct or inverse, i.e.
there exists a degree function from the set of objects of
$\mathcal{B}$ to some ordinal such that every non-identity map of
$\mathcal{B}$ always raises or always lowers the degree. Then the
canonical functor
\[\dgm_\mathcal{B} : \ho(\mathcal{M}^\mathcal{B}) \longrightarrow
\ho(\mathcal{M})^\mathcal{B}\] from the homotopy category of diagrams
of objects of $\mathcal{M}$ over $\mathcal{B}$ to the category of
diagrams of objects of $\ho(\mathcal{M})$ over $\mathcal{B}$ is full
and essentially surjective.  \eth

The homotopy category of flows $\ho(\dtop)$ is weakly complete and
weakly cocomplete as any homotopy category of any model category
\cite{MR99h:55031}. Weak limit and weak colimit satisfy the same
property as limit and colimit except the uniqueness.  Weak small
(co)products coincide with small (co)products. Weak (co)limits can be
constructed using small (co)products and weak (co)equalizers in the
same way as (co)limits are constructed by small (co)products and
(co)equalizers (\cite{MR1712872} Theorem~1 p109). And a weak
coequalizer
\[A\stackrel{f,g}\rightrightarrows B \stackrel{h}\longrightarrow D\] 
is given by a weak pushout 
\[
\xymatrix{
B \fr{h} && D \\
&& \\
A\sqcup B \fu{(f,\id_B)}\fr{(g,\id_B)} && B. \fu{h}}
\] 
And finally, weak pushouts (resp. weak pullbacks) are given by
homotopy pushouts (resp. homotopy pullbacks) (e.g., \cite{rep}
Remark~4.1 and \cite{MR920963} Chapter~III). As explain in
\cite{cofibrationcat}, Theorem~\ref{weak} can be also used for the
construction of certain kind of weak limits and of weak colimits:

\begin{cor} \label{weakweak} (\cite{cofibrationcat} Theorem~8.8.6) Let
  $\mathcal{M}$ be a model category. Let $\mathcal{B}$ be a small
  Reedy category which is free, i.e. freely generated by a graph.
  Moreover, let us suppose that $\mathcal{B}$ is either direct or
  inverse. Let $X\in \ho(\mathcal{M})^\mathcal{B}$. Let $X'\in
  \mathcal{M}^\mathcal{B}$ with $\dgm_\mathcal{B}(X')=X$.
\begin{enumerate}
\item If $\mathcal{B}$ is direct, then a weak colimit $\wliminj X$
  of $X$ is given by \[\boxed{\wliminj X := \dgm_\mathcal{B}(\holiminj X')
  \simeq \dgm_\mathcal{B}(\liminj X'^{cof})}\] where the cofibrant
  replacement $X'^{cof}$ is taken in the Reedy model structure of
  $\mathcal{M}^\mathcal{B}$. This weak colimit, called the privileged
  weak colimit in Heller's terminology, is unique up to a
  non-canonical isomorphism.
\item If $\mathcal{B}$ is inverse, then a weak limit $\wlimproj X$ of
  $X$ is given by \[\boxed{\wlimproj X := \dgm_\mathcal{B}(\holimproj
    X') \simeq \dgm_\mathcal{B}(\limproj X'^{fib})}\] where the
  fibrant replacement $X'^{fib}$ is taken in the Reedy model structure
  of $\mathcal{M}^\mathcal{B}$. This weak limit, called the privileged
  weak limit in Heller's terminology, is unique up to a non-canonical
  isomorphism.
\end{enumerate}
\end{cor}

We have now the necessary tools to state the theorem:

\bth \label{weakweakweak} For each process name $P$ of our restriction
of CCS, consider the object $[P]$ of $\ho(\dtop)$ defined by induction
on the syntax of $P$ as in Table~\ref{wsem}. Then one has $\llbracket
P\rrbracket \in [P]$, i.e. the weak S-homotopy type of $\llbracket
P\rrbracket$ is $[P]$. \eth

\bpf One observes that the small categories involved for the
construction of pushouts and colimits of towers are Reedy direct free
and that the small category involved for the construction of pullbacks
is Reedy inverse free. One then proves $\llbracket P\rrbracket \in
[P]$ by induction on the syntax of $P$ with Corollary~\ref{weakweak},
Proposition~\ref{L1}, Proposition~\ref{L2} and Corollary~\ref{L3}.
\epf

We do not know how to construct a pure homotopical semantics of the
parallel composition with synchronization.

\end{document}